

\magnification=\magstep1
\baselineskip=18 true pt


\def\COMPARE#1#2#3#4{{\edef\worda{#1}\edef\wordb{#2}%
\ifx\worda\wordb#3\else#4\fi}}

\def\rmcolon[#1:]{#1}
\def\rmslash[#1/]{#1}
\def\takeone[#1/#2]{\COMPARE{#2}{}{#1}{\rmslash[#2]}\COMPARE{#2}{/}{/#1/}{}}
\def\refref[#1:#2]{\COMPARE{#2}{}{{\rm (\takeone[#1/])}}{\rmcolon[#2]}}

\def\label[#1]{\ref[#1]}
\def\labeldef[#1]{\ref[#1]}
\def\ref[#1]{\refref[#1:]}

\def\raw[#1]{#1}

\def\citecite[#1:#2]{\COMPARE{#2}{}{{\bf #1}}{{\bf \rmcolon[#2]}}}

\def\bibitem[#1]{\item{[\cite[#1]]}}
\def\cite[#1]{\citecite[#1:]}

\def\citeinitial{\def\citecite[##1:##2]{{\bf ##1}}}

\citeinitial


\def\BScap#1#2{\bigbreak\bookmark[]{1}{#2}\noindent{\bf #1}\nobreak\smallskip\noindent\ignorespaces}%
\def\BS#1{\BScap{#1}{#1}}%

\def\BSScap#1#2{%
\medbreak\bookmark[]{2}{#2}\noindent{\bf #1}\nobreak\smallskip\noindent\ignorespaces}  

\def\BSS#1{\BSScap{#1}{#1}}

\def\bibliography{\bookmark[]{1}{References}\centerline{\bf References}}


\def\claim#1. {\noindent{\bf #1. }}

\def\flushright#1{{\unskip\nobreak\hfil\penalty50\hskip2em\hbox{}\nobreak\hfil%
#1\parfillskip=0pt\finalhyphendemerits=0\par}}

\def\bull{\vrule height 1.8ex width 1.0ex depth .1ex }
\def\QED{\ifmmode\eqno\hbox{$\bull$}\else\flushright{\hbox{$\bull$}}\fi}


\def\R{\mathop{\bf R}\nolimits}
\def\N{\mathop{\bf N}\nolimits}

\def\abs#1{|{#1}|}

\def\norm#1{\|{#1}\|}

\def\dist{\mathop{\rm dist}\,}

\def\half{{1 \over 2}}


\def\m{\hfil\break}  




\def\bkmk[#1]#2#3{}%


\def\bookmark[#1]#2#3{{\def\label[##1:##2]{##2}\def\labeldef[##1:##2]{##2}%
\def\ref[##1:##2]{##2}\bkmk[#1]{#2}{#3}}}

\ifx\pdftexversion\pdftexrevision\else\errmessage{Please use plain tex not pdftex}\fi


\def\epsilon{\varepsilon}
\def\calE{{\cal E}}


\null
\vskip 0.8 true cm

\centerline{\bf Remarks on the Clark theorem}

\bigskip

\centerline{Guosheng Jiang${}^{1}$, Kazunaga Tanaka${}^{2}$, Chengxiang Zhang${}^{3}$}

\bigskip

\settabs 25 \columns
\+&$1$ & School of Mathematical Sciences, Capital Normal University\cr
\+&&Beijing 100048, China\cr
\+&$2$ & Department of Mathematics, School of Science and Engineering\cr
\+&& Waseda University, 3-4-1 Ohkubo, Shinjuku-ku, Tokyo 169-8555, Japan\cr
\+&$3$ & Chern Institute of Mathematics and LPMC, Nankai University\cr
\+&&Tianjin 300071, China \cr

\bigskip
\bigskip

{\narrower
\noindent{\bf Abstract.}
The Clark theorem is important in critical point theory.  For a class of even functionals it ensures 
the existence of infinitely many negative critical values converging to $0$ and 
it has important applications to sublinear elliptic problems.   
We study the convergence of the corresponding critical points and we give a characterization of 
accumulation points of critical points together with examples, in which critical points with negative 
critical values converges to non-zero critical point.  Our results improve the abstract results 
in Kajikiya [\cite[Ka1]] and Liu-Wang [\cite[LW]].

}

\bigskip
\bigskip

\BS{1. Introduction and main results}
The Clark theorem is one of the most important results in critical point theory 
(Clark [\cite[Cl]], see also Heinz [\cite[H]]).  It was successfully applied to 
sublinear elliptic problems with odd symmetry and the existence
of infinitely many solutions which accumulate to $0$ was shown.

To state the Clark theorem, we need some terminologies:
let $(X,\norm\cdot_X)$ be a Banach space and $I\in C^1(X,\R)$.

\smallskip

\item{(i)} For $c\in\R$ we say that $I(u)$ satisfies the $(PS)_c$ condition if any sequence $(u_j)_{j=1}^\infty\subset X$
with $I(u_j)\to c$, $\norm{I'(u_j)}_{X^*}\to 0$ has a convergent subsequence.
\item{(ii)} Let ${\cal E}$ be the family of sets $A\subset X\setminus\{ 0\}$ such that $A$ is closed and
symmetric with respect to $0$.
For $A\in\calE$, the genus $\gamma(A)$ is introduced by Krasnosel'skii [\cite[Kr]]
(c.f. Coffman [\cite[Co]], Rabinowitz [\cite[R]]) as the smallest integer $n$ such that there exists 
an odd continuous map $\zeta\in C(A,\R^n\setminus\{ 0\})$.  When there does not exist such a map, we set
$\gamma(A)=\infty$.  See Rabinowitz [\cite[R]] for fundamental properties of the genus.

\smallskip

\noindent
Now we give a variant of the Clark theorem due to Heinz [\cite[H]].

\medskip

\proclaim Theorem 1.1 {\rm (Heinz [\cite[H]])}.  
Let $(X,\norm\cdot_X)$ be a Banach space and suppose that $I(u)\in C^1(X,\R)$ satisfies the following conditions:
{\parindent=1.5\parindent
\item{\rm (A1)} $I(0)=0$.  $I(u)$ is even in $u$ and bounded from below;
\item{\rm (A2)} $I(u)$ satisfies $(PS)_c$ for all $c<0$;
\item{\rm (A3)} For any $k\in\N$, there exists $A\in{\cal E}$ such that
    $$  \gamma(A)\geq k \quad \hbox{and}\quad \sup_{u\in A}I(u)<0.
    $$
}
Then $I(u)$ has a sequence $(c_j)_{j=1}^\infty$ of critical values of $I(u)$ such that
    $$  \eqalign{
        &c_j<0 \quad \hbox{for all}\ j\in \N,\cr
        &c_j\to 0 \quad \hbox{as}\ j\to\infty.\cr}
    $$
Here 
    $$  c_j=\inf_{A\in{\cal E}, \gamma(A)\geq j}\sup_{u\in A} I(u).         \eqno\label[1.1]
    $$

\medskip

\claim Remark 1.2.
In [\cite[H]], it was assumed that 

{\parindent=1.5\parindent
\item{\rm (A2')} $I(u)$ satisfies $(PS)_c$ for all $c\in\R$.
}

\noindent
From its proof, we can easily see that $(PS)_c$ just for $c<0$ is enough for the existence of critical values.

\medskip

By Theorem 1.1, there exists a sequence $(u_j)_{j=1}^\infty$ of critical points of $I(u)$ such that
$I(u_j)=c_j\to -0$ as $j\to\infty$.  Thus it is natural to ask whether $u_j\to 0$ holds or not.
More generally, the existence of a sequence of non-zero critical points $(u_j)_{j=1}^\infty$ (or critical points
with negative critical values) satisfying $u_j\to 0$ is of interest.
This question has been studied  by Kajikiya [\cite[Ka1]] and Liu-Wang [\cite[LW]] together with applications
to sublinear elliptic problems.  We note that Liu-Wang [\cite[LW]] also studied periodic solutions of
Hamiltonian systems.  More precisely, under the assumptions of (A1), (A2') and (A3),  Kajikiya [\cite[Ka1]]
showed either
{\parindent=1.5\parindent
\item{(C1)} There exists a sequence $(u_j)_{j=1}^\infty$ such that
    $$  I'(u_j)=0, \ I(u_j)<0 \ \hbox{and}\ u_j\to 0 \ \hbox{as}\ j\to \infty.
    $$
}

\noindent
or
{\parindent=1.5\parindent
\item{(C2)} There exists two sequences $(u_j)_{j=1}^\infty$ and $(v_j)_{j=1}^\infty$ such that
    $$  I'(u_j)=0, \ I(u_j)=0, \ u_j\not=0 \ \hbox{and}\ u_j\to 0 \ \hbox{as}\ j\to \infty
    $$
and
    $$  I'(v_j)=0, \ I(v_j)<0 \ \hbox{and $v_j$ converges to a non-zero limit.}
    $$
}

\noindent
holds.

Liu-Wang [\cite[LW]] assumed (A1), (A2') and the following (A3'), which is stronger than (A3), 

{\parindent=1.5\parindent
\item{\rm (A3')} For any $k\in \N$ there exists a $k$-dimensional subspace $X^k$ of $X$ and $\rho_k>0$ such that
    $$  \sup\{ I(u);\, u\in X^k,\ \norm u_X =\rho_k\} <0
    $$
}

\noindent
and they showed either (C1) above or

{\parindent=1.5\parindent
\item{(C3)} There exists $r>0$ such that for any $0<a<r$ there exists a critical point $u$ such that
    $$  \norm u_X=a \quad \hbox{and}\quad I(u)=0.
    $$
}

\medskip

\noindent
In what follows, we denote by $\widehat K_0$ the connected component of 
$K_0=\{ u\in X;\, I'(u)=0, \, I(u)=0\}$ including $0$.

\medskip

\claim Remark 1.3.
From their proof of their main result, Liu-Wang [\cite[LW]] claimed that (C3) can be strengthened as 

{\parindent=1.5\parindent

\item{(C3')} There exists $r>0$ such that 
    $$  \widehat K_0 \cap \{ u\in X;\, \norm u_X=r\} \not= \emptyset.
    $$

}

\medskip

The aim of this paper is to show the following Theorem 1.4 and Theorem 1.6;  In Theorem 1.4, 
we give a new characterization of accumulation points of critical points with negative critical values 
and unifies the results of Kajikiya and Liu-Wang.  
On the other hand, in Theorem 1.6 we answer a natural question concerning (C1), which is stated below.
We believe that Theorems 1.4 and 1.6 give us a better understanding of the Clark theorem.

First we give our Theorem 1.4.

\medskip

\proclaim Theorem 1.4.  
Let $(X,\norm\cdot_X)$ be a Banach space and suppose $I\in C^1(X,\R)$ satisfies {\rm (A1)}, {\rm (A3)} and
{\parindent=1.5\parindent
\item{\rm (A2'')} $I(u)$ satisfies $(PS)_c$ for all $c\leq 0$.
}

\noindent
{\sl Then there exists a sequence $(u_j)_{j=1}^\infty\subset X$ of critical points of $I(u)$ such
that
    $$  \eqalignno{
        &I(u_j)<0 \quad \hbox{for all}\ j\in\N,         &\label[1.2]\cr
        &I(u_j) \to 0 \quad  \hbox{as}\ j\to\infty,     &\label[1.3]\cr
        &\dist(u_j,\widehat K_0)\,\left(\equiv \inf\{ \norm{u_j-v}_X;\, v\in \widehat K_0\}\right)
            \to 0 \quad \hbox{as}\ j\to\infty.\cr}
    $$
}
\medskip

\noindent
As an immediate corollary to our Theorem 1.4, we have

\medskip

\proclaim Corollary 1.5.
Under the assumptions of Theorem 1.4, assume that (C1) does not take  place.
Then $\widehat K_0\not= \{ 0\}$.

\medskip

\noindent
Since $\widehat K_0\not=\{ 0\}$ implies (C2) and (C3), Corollary 1.5 covers the results of 
Kajikiya [\cite[Ka1]] and Liu-Wang [\cite[LW]].

\medskip

Next we study a question concerning (C1).  
In many applications of the Clark theorem to sublinear elliptic equations, there exist sequences 
$(u_j)_{j=1}^\infty$ of solutions with  \ref[1.2], \ref[1.3] and
    $$  u_j\to 0 \quad \hbox{as}\ j\to\infty.                \eqno\label[1.4]
    $$
So (C1) may be expected under the assumption of Theorem 1.4 and a natural question is to ask 
whether (C1) always takes place under the assumption of Theorem 1.4 or not.
Our Theorem 1.6 answers this question negatively.

\medskip

\proclaim Theorem 1.6. 
Conditions (A1), (A2''), (A3') do not imply (C1).  In particular, under the assumptions of Theorem 1.4, 
(C1) does not hold in general.

\medskip

\claim Remark 1.7.
An example related to our Theorem 1.6 was given in Example \raw[1.3] of [\cite[Ka1]] (c.f. [\cite[Ka2]]).
It shows that there exists a functional $I\in C^1(X,\R)$ which satisfies (A1), (A2''), (A3) and the 
following property:
\itemitem{} There exists an $r_0>0$ independent of $j$ such that
    $$  I'(u)=0 \ \hbox{and}\ I(u)=c_j \quad \hbox{imply} \quad \norm u_X\geq r_0.
    $$

\noindent
Here $c_j$ is given in \ref[1.1] and $c_j$ satisfies $c_j<0$ and $c_j\to 0$ as 
$j\to \infty$.  Thus a special case of (C1) does not hold for $I$.  
In Section \ref[Subsection:3.1] we give another example $I\in C^1(\ell^2,\R)$ for which 
we give an explicit description of all critical points of $I(u)$ and no critical points with negative 
critical values do not exist in a neighborhood of $0$.  Especially (C1) does not hold for our $I(u)$.
Our example also shows a typical situation of our Theorem 1.4.

\medskip

Finally we remark that in our Theorem 1.4, (A2''), especially $(PS)_0$ is important.  In fact, we have

\medskip

\proclaim Theorem 1.8.
Under the assumptions of Theorem 1.1, especially without $(PS)_0$, the conclusion of Theorem 1.4 does not
hold in general.

\medskip

In the following Section 2, we give a proof to our Theorem 1.4.  Here estimates of $I'(u)$ play important
roles.  In Section 3, we give two examples which show Theorems 1.6 and 1.8.

\medskip

\BS{2. Proof of Theorem 1.4}
In what follows, we use the following notation for $\delta>0$
    $$  \eqalign{
        B_\delta(u)&=\{ x\in X;\, \norm{x-u}_X<\delta\} \quad\hbox{for}\ u\in X,\cr
        N_\delta(D)&=\{ x\in X;\, \dist(x,D)<\delta\} 
            \quad \hbox{for}\ D\subset X,\cr}
    $$
where
    $$  \dist(x,D)=\inf_{y\in D}\norm{x-y}_X.
    $$
We note that $N_\delta(D)=\bigcup_{y\in D} B_\delta(y)$.

\BSS{2.1. A fundamental fact from topology}
To show our Theorem 1.4, we need the following characterization of connected components of
compact sets.

\proclaim Lemma 2.1.
Let $D\subset X$ be a compact set such that $0\in D$.  For $\delta>0$, let $O_\delta$
be the connected component of $N_\delta(D)$ including $0$.  Then we have
    $$  \bigcap_{\delta>0} \overline{O_\delta}=\widehat D,
    $$
where $\widehat D$ is the connected component of $D$ including 0.

\claim Proof.
By the definition of $O_\delta$ and $\widehat D$, it is clear that $\widehat D\subset O_\delta$
for all $\delta>0$.  Thus
    $$  \widehat D\subset \bigcap_{\delta>0} O_\delta
        \subset \bigcap_{\delta>0} \overline{O_\delta}.
    $$
By the compactness of $D$, we also have $D=\bigcap_{\delta>0} \overline{N_\delta(D)}$.

We set 
    $$  A=\bigcap_{\delta>0} \overline{O_\delta}\subset D.
    $$
It suffices to show that $A$ is connected.  For $\delta>0$ we also set 
$D_\delta=\overline{O_\delta}\cap D$.  Then we have
    $$  \eqalignno{
        &A = \bigcap_{\delta>0} D_\delta,                   &\label[2.1]\cr
        &O_\delta = \bigcup_{u\in D_\delta} B_\delta(u),    &\label[2.2]\cr
        &\hbox{$\delta_1<\delta_2$ implies $D_{\delta_1}\subset D_{\delta_2}$.} &\label[2.3]\cr}
    $$
Arguing indirectly, we suppose that $A$ is not connected.  Then there exist two compact
sets $A_1$, $A_2\subset X$ such that $A_1\cap A_2=\emptyset$, $A_1\cup A_2=A$.  We set
    $$  \beta=\half\dist(A_1,A_2)>0.
    $$
For each $\delta>0$, since $O_\delta$ is a connected set including $A_1\cup A_2$, we
have
    $$  O_\delta\cap\{ x\in X;\, \dist(x,A_1)=\beta\} \not=\emptyset.
    $$
By \ref[2.2], we can see that for any $x\in O_\delta$ there exists $u\in D_\delta$ such that
$x\in B_\delta(u)$.  Thus
    $$  D_\delta\cap\{ x\in X;\, \dist(x,A_1)\in [\beta-\delta,\beta+\delta]\} \not=\emptyset.
    $$
Since $\Bigl( D_\delta\cap\{ x\in X;\, \dist(x,A_1)\in [\beta-\delta,\beta+\delta]\} \Bigr)_{\delta>0}$
has the finite intersection property by \ref[2.3], we have
    $$  \eqalign{
        & A\cap \{ x\in X;\, \dist(x,A_1)=\beta\}\cr
        =& \bigcap_{\delta>0} \Bigl( D_\delta\cap\{ x\in X;\, \dist(x,A_1)\in [\beta-\delta,\beta+\delta]\} \Bigr)
            \not=\emptyset,\cr}
    $$
which contradicts with the choice of $\beta>0$.  Thus $A$ is a connected set.  \QED

\medskip

\BSS{2.2. A gradient estimate}
Suppose that $I(u)\in C^1(X,\R)$ satisfies the assumptions of Theorem 1.4.  We use the
following notation:
    $$  \eqalign{
        K&=\{ u\in X;\, I'(u)=0\},\cr
        K_0 &= \{ u\in K;\, I(u)=0\}, \cr
        K_- &= \{ u\in X;\, \hbox{there exists $(v_j)_{j=1}^\infty\subset K$ such that}\cr
        & \qquad\qquad  \ \ I(v_j)<0 \ \hbox{for all $j$ and $I(v_j)\to 0$, $v_j\to u$ as $j\to\infty$} \}.\cr}
    $$
By $(PS)_0$, we have $K_-\subset K_0$.  We also use notation for $a<b$
    $$  \eqalign{
        &[I\leq a] =\{ u\in X;\, I(u)\leq a\},\cr
        &[a\leq I \leq b] =\{ u\in X;\, a\leq I(u)\leq b\}.\cr}
    $$
It is clear that $0\in K_0$.  We denote  by $\widehat K_0$ the connected component of $K_0$
including $0$.  To show our Theorem 1.4 it suffices to prove
    $$  K_-\cap \widehat K_0 \not=\emptyset.                    \eqno\label[2.4]
    $$
For $\delta>0$, let $O_\delta$ be the connected component of $N_\delta(K_0)$ including $0$.
By Lemma 2.2, we have
    $$  \widehat K_0 = \bigcap_{\delta>0}\overline{O_\delta}.
    $$
Thus to prove \ref[2.4] it suffices to show
    $$  \overline{O_\delta} \cap K_- \not=\emptyset \quad \hbox{for all}\ \delta>0.
    $$
We argue indirectly and suppose for some $\delta_0>0$
    $$  \overline{O_{\delta_0}} \cap K_- =\emptyset.        \eqno\label[2.5]
    $$
Under the assumption \ref[2.5], we set
    $$  K_{0,i}= \overline{O_{\delta_0}}\cap K_0, \quad
        K_{0,e}= K_0\setminus O_{\delta_0}.
    $$
Then $K_{0,i}$ and $K_{0,e}$ are disjoint compact sets such that $K_0=K_{0,i}\cup K_{0,e}$ and
    $$  \eqalignno{
        &\dist(K_{0,i}, K_{0,e}) \geq 2\delta_0,            &\label[2.6]\cr
        &K_- \subset K_{0,e}.                               &\label[2.7]\cr}
    $$
We note that \ref[2.7] follows from \ref[2.5].

First we have

\proclaim Lemma 2.2.
Assume \ref[2.5].  Then for any $r>0$ there exist $\rho>0$ and $\nu>0$ such that
    $$  \eqalignno{
        &[-\rho\leq I\leq 0] \cap K\subset N_r(K_0),            &\label[2.8]\cr
        &[-\rho\leq I<0] \cap K \subset N_r(K_-)\subset N_r(K_{0,e}),   &\label[2.9]\cr
        &\norm{I'(u)}_{X^*}\geq \nu \quad 
            \hbox{for all}\ u\in [-\rho\leq I\leq 0]\setminus N_r(K_0). &\label[2.10]\cr}
    $$
Moreover for any $\epsilon\in (0,\rho)$ there exists $\nu_\epsilon\in (0,\nu]$ such that
    $$  \norm{I'(u)}_{X^*} \geq \nu_\epsilon \quad \hbox{for  all}\ 
        u\in [-\rho\leq I\leq -\epsilon]\setminus N_r(K_{0,e}). \eqno\label[2.11]
    $$

\claim Proof.
Using $(PS)_0$ and the definition of $K_-$, we can check \ref[2.8]--\ref[2.10] easily
for small $\rho$ and $\nu>0$.  
We show \ref[2.11].  Suppose that for $r$, $\rho$, $\nu>0$, \ref[2.8]--\ref[2.10] hold.
If \ref[2.11] does not hold, we can find $\epsilon \in (0,\rho)$ and a sequence
$(u_j)_{j=1}^\infty$ such that
    $$  I(u_j) \in [-\rho,-\epsilon],\quad \norm{I'(u_j)}\to 0, \quad
        u_j\not\in N_r(K_{0,e}),
    $$
By $(PS)$, we can extract a subsequence $(u_{j_k})$ such that $u_{j_k}\to u_0$
for some $u_0\in [-\rho\leq I\leq -\epsilon]\cap K$.  By \ref[2.9], we have
$u_{j_k}\not\in N_r(K_{0,e})$ for large $k$, which is a contradiction.  
Thus we have \ref[2.11].
\QED

\medskip

\BSS{2.3.  Deformation argument}
The aim of this section is the following

\proclaim Proposition 2.3.  
Assume \ref[2.5].  Then for any $r\in (0,\delta_0/3]$ there exists $d>0$ with the following
property: for any $\epsilon\in (0,d/2]$ there exists an odd continuous map $\eta_\epsilon:\,
[I<0]\to [I<0]$ such that
    $$  \eta_\epsilon([I\leq-\epsilon]) \subset [I\leq -d] \cup N_{3r}(K_{0,e}).
    $$

\claim Proof.
First we define an ODE in $X$ to define $\eta_\epsilon$. \m
For a given $r\in (0,\delta_0/3]$, let $\rho$, $\nu>0$ be constants given in Lemma 2.2.  We set
    $$  d = {1\over 3}\min\{ \rho, \nu r\}>0.               \eqno\label[2.12]
    $$
Then again by Lemma 2.2, for any given $\epsilon\in (0,d]$ there exists $\nu_\epsilon>0$ with 
the property \ref[2.11].

By \ref[2.9], we have $I'(u)\not=0$ for all $u\in [-3d\leq I<0]\setminus N_r(K_{0,e})$.  
Thus there exists a locally Lipschitz odd vector field $V(u):\, [-3d\leq I<0]\setminus N_r(K_{0,e})
\to X$ such that
    $$  \eqalignno{
        &\norm{V(x)}_X \leq 1 \quad \hbox{for all}\ u\in [-3d\leq I<0]\setminus N_r(K_{0,e}), \cr
        &I'(u)V(u) >0 \quad \hbox{for all}\ u\in [-3d\leq I<0]\setminus N_r(K_{0,e}), &\label[2.13]\cr
        &I'(u)V(u) \geq {\nu\over 2} \quad \hbox{for all}\ u\in [-3d\leq I<0]\setminus N_r(K_0), &\label[2.14]\cr
        &I'(u)V(u) \geq {\nu_\epsilon\over 2} \quad \hbox{for all}\ u\in [-3d\leq I\leq -\epsilon]\setminus N_r(K_{0,e}).
                                                        &\label[2.15]\cr}
    $$
Let $\phi_1(u)$, $\phi_2(u):\, X\to [0,1]$ be even Lipschitz continuous functions such that
    $$  \phi_1(u)=\cases{   1   &for $u\in [-d\leq I]$,\cr
                            0   &for $u\in [I\leq -2d]$,\cr} \qquad
        \phi_2(u)=\cases{   1   &for $u\in X\setminus N_{2r}(K_{0,e})$,\cr
                            0   &for $u\in N_r(K_{0,e})$.\cr}
    $$
We set
    $$  \widetilde V(u)=\phi_1(u)\phi_2(u)V(u)
    $$
and we note that $\widetilde V(u)$ is well-defined on $[I<0]$.  For $u\in [I<0]$ we consider
    $$  \left\{\eqalign{
        &{d\eta\over dt} = -\widetilde V(\eta), \cr
        &\eta(0,u)=u. \cr}
        \right. 
    $$
We have for all $(t,u)$
    $$  \eqalignno{
        &\eta(t,-u)=-\eta(t,u),                 &\label[2.16]\cr
        &\norm{{d\eta\over dt}}_X \leq 1.       &\label[2.17]\cr}
    $$
Since 
    $$  {d\over dt}I(\eta(t,u)) = I'(\eta(t,u)){d\eta\over dt} = -I'(\eta)\widetilde V(\eta),
    $$
it follows from \ref[2.13]--\ref[2.15] that
    $$  \eqalignno{
        &{d\over dt}I(\eta(t,u)) \leq 0 \quad \hbox{if}\ \eta(t,u)\in [I<0],    &\label[2.18]\cr
        &{d\over dt}I(\eta(t,u)) \leq -{\nu\over 2} \quad \hbox{if}\ 
            \eta(t,u)\in [-d\leq I<0]\setminus N_{2r}(K_0), &\label[2.19]\cr
        &{d\over dt}I(\eta(t,u)) \leq -{\nu_\epsilon\over 2} \quad \hbox{if}\ 
            \eta(t,u)\in [-d\leq I\leq -\epsilon]\setminus N_{2r}(K_{0,e}). &\label[2.20]\cr}
    $$
By \ref[2.17] and \ref[2.18], we note that for any $u\in [I<0]$, $\eta(t,u)$ exists globally, that is,
$\eta(t,u):\, [0,\infty)\times [I<0]\to [I<0]$ is well-defined.  
For a latter use, we note that
    $$  \overline{N_{3r}(K_{0,e})}\setminus N_{2r}(K_{0,e}) \subset X\setminus N_{2r}(K_{0}).
    $$
Thus by \ref[2.19]
    $$  {d\over dt}I(\eta(t,u)) \leq -{\nu\over 2} \quad \hbox{if}\ 
            \eta(t,u)\in [-d\leq I<0]\cap (N_{3r}(K_{0,e})\setminus N_{2r}(K_{0,e})).
                                                        \eqno\label[2.21]
    $$
Next we claim that

\medskip

\claim Claim.  Let $T_\epsilon={2d\over\nu_\epsilon}$.  Then 
    $$  \eta(T_\epsilon,u)\in [I\leq -d]\cup N_{3r}(K_{0,e}) \quad
                \hbox{for any}\  u\in [I\leq -\epsilon].            \eqno\label[2.22]
    $$

\noindent
To prove \ref[2.22], it suffices to show that if $u\in [I\leq -\epsilon]$ satisfies
    $$  \eta(T_\epsilon,u)\not\in [I\leq -d],       \eqno\label[2.23]
    $$
then
    $$  \eta(T_\epsilon,u)\in N_{3r}(K_{0,e}).      \eqno\label[2.24]
    $$
We note that under the condition \ref[2.23]
    $$  \eta(t,u) \in [-d\leq I\leq -\epsilon] 
        \quad \hbox{for all}\ t\in [0,T_\epsilon]. \eqno\label[2.25]
    $$

\smallskip

\noindent
{\sl Step 1: Assume \ref[2.23], i.e.,  \ref[2.25].  Then
    $$  \eta([0,T_\epsilon],u) \cap N_{2r}(K_{0,e}) \not= \emptyset.    \eqno\label[2.26]
    $$
}

\smallskip

\noindent
In fact, if \ref[2.26] does not hold, it follows from \ref[2.20] that
    $$  {d\over ds}I(\eta(s,u)) \leq -{\nu_\epsilon\over 2} \quad \hbox{for all}\ s\in [0,T_\epsilon].
    $$
Thus, by the definition of $T_\epsilon$,
    $$  \eqalign{
        I(\eta(T_\epsilon,u)) &\leq I(u) + \int_0^{T_\epsilon} {d\over ds}I(\eta(s,u))\, ds\cr
            & \leq -\epsilon -{\nu_\epsilon\over 2} T_\epsilon \cr
            & < -d, \cr}
    $$
which is in contradiction with \ref[2.23].

\smallskip

\noindent
{\sl Step 2:  Assume \ref[2.23], i.e., \ref[2.25].  Then \ref[2.24] holds.}

\smallskip

\noindent
Assume  \ref[2.24] does not hold.  Then $\eta(T_\epsilon,u) \not\in N_{3r}(K_{0,e})$ and by
\ref[2.26] the orbit $\eta(t,u)$ enters in $N_{2r}(K_{0,e})$ for some $t\in [0,T_\epsilon]$.
Thus there exists an interval $[t_0,t_1] \subset [0,T_\epsilon]$ such that
    $$  \eqalign{
        &\eta(t_0,u)\in \partial N_{2r}(K_{0,e}), \cr
        &\eta(t_1,u)\in \partial N_{3r}(K_{0,e}), \cr
        &\eta(t,u)\in \overline{N_{3r}(K_{0,e})}\setminus N_{2r}(K_{0,e}) \quad
            \hbox{for all}\ t\in [t_0,t_1]. \cr}
    $$
By \ref[2.17], we have
    $$  \eqalign{
        r &\leq \norm{\eta(t_1,u)-\eta(t_0,u)}_X \leq \int_{t_0}^{t_1} \norm{{d\over ds}\eta(s,u)}_X\, ds\cr
            &\leq t_1-t_0.\cr}
    $$
Thus by \ref[2.21],
    $$  \eqalign{
        I(\eta(T_\epsilon,u)) &\leq I(u) + \int_0^{T_\epsilon} {d\over ds}I(\eta(s,u))\,ds\cr
        &\leq I(u) + \int_{t_0}^{t_1} {d\over ds}I(\eta(s,u))\, ds\cr
        &\leq I(u) -{\nu\over 2} (t_1-t_0)\cr
        &\leq -\epsilon -{\nu r\over 2} \leq -d.\cr}
    $$
This is a contradiction to \ref[2.23].  Thus we have $\eta(T_\epsilon,u)\in N_{3r}(K_{0,e})$ and
the conclusion of Step 2 holds.

Setting $\eta_\epsilon(u)=\eta(T_\epsilon,u)$, we have the desired deformation.  \QED

\medskip

\BSS{2.4. End of the proof of Theorem 1.4}

\claim Proof of Theorem 1.4.
Since $K_{0,e}\in\calE$ is compact, we can see
    $$  \gamma_{0,e}\equiv \gamma(K_{0,e})<\infty.
    $$
Moreover for small $r\in (0,\delta_0/3]$ 
    $$  \gamma(\overline{N_{3r}(K_{0,e})}) = \gamma(K_{0,e}) = \gamma_{0,e}.
    $$
We fix such an $r$ and we choose $d>0$ by Proposition 2.3.  

By Clark's theorem [\cite[Cl]], we have
    $$  c_k\equiv \inf_{\gamma(A)\geq k}\sup_{u\in A} I(u) \nearrow 0.
    $$
Thus there exists $k_0$ such that $c_{k_0} > -d$.  That is, $\gamma([I\leq -d]) < k_0$.
Thus
    $$  \gamma([I\leq -d]\cup \overline{N_{3r}(K_{0,e})})
        \leq \gamma([I\leq -d]) +\gamma( \overline{N_{3r}(K_{0,e})}) 
        \leq \gamma_{0,e}+ k_0.                 \eqno\label[2.27]           
    $$
By the assumption of Theorem 1.4, there exists $A\in\calE$ such that
    $$  \gamma(A) > \gamma_{0,e} + k_0 \quad \hbox{and}\quad
        \sup_{u\in A} I(u) <0.  
    $$
Choosing $\epsilon\in (0,d)$ such that $\sup_{u\in A}I(u) <-\epsilon$, we have
    $$  \gamma([I\leq -\epsilon]) \geq \gamma(A) > \gamma_{0,e}+ k_0.   \eqno\label[2.28]
    $$
On the other hand, by Proposition 2.3, there exists a continuous odd map $\eta_\epsilon:\,
[I<0]\to [I<0]$ such that
    $$  \eta_\epsilon([I\leq -\epsilon])\subset [I\leq -d]\cup N_{3r}(K_{0,e}).
    $$
Thus by \ref[2.27]
    $$  \eqalign{
        \gamma([I\leq -\epsilon]) &\leq \gamma(\eta_\epsilon([I\leq -\epsilon])) \cr
        &\leq \gamma([I\leq -d]\cup \overline{N_{3r}(K_{0,e})}) \cr
        &\leq \gamma_{0,e}+k_0,\cr}
    $$
which is in contradiction with \ref[2.28].  Thus \ref[2.5] cannot take place and we 
complete the proof of our Theorem 1.4.  \QED

\medskip

\BS{3. Some examples}
In this section we give two examples which show Theorems 1.6 and 1.8.

\medskip

\BSS{3.1. An example which shows Theorem 1.6}
We give an example which shows that (A1), (A2''), (A3') do not imply (C1).  
We work in the space $X=\ell^2$, that is,
    $$  \eqalign{
        &X = \left\{ (t,x_1,x_2,\cdots);\, t,\, x_j \in\R \ (j=1,2,\cdots), \ 
            t^2+\sum_{j=1}^\infty \abs{x_j}^2 <\infty\right\},\cr
        &\norm{(t,x_1,x_2,\cdots)}_X = \left(t^2+\sum_{j=1}^\infty \abs{x_j}^2\right)^{1/2}.\cr}
    $$
Since the first component has a special role in our argument, we use notation 
$(t,x_1,x_2,\cdots)$ for elements of $X$.

We consider a functional $I(t,x_1,x_2,\cdots):\, X\to \R$ in a form
    $$  I(t,x_1,x_2,\cdots)= \half \sum_{j=1}^\infty \abs{x_j}^2 
        -{2\over 3}\sum_{j=1}^\infty 3^{-j} \left(a_+(t)(x_j)_+^{3/2}+a_-(t)(x_j)_-^{3/2}\right) +\varphi(t),
    $$
where $x_+=\max\{ x, 0\}$, $x_-=\max\{ -x, 0\}$ and $a_+(t)$, $a_-(t)\in C^1(\R,\R)$ are given by
    $$  a_+(t) =2 +\mu(t), \quad a_-(t)=2-\mu(t).
    $$
Here $\mu(t)\in C^1(\R,\R)$ satisfies
    $$  \eqalignno{
        \mu(t) &\in [-1,1] \quad \hbox{for all}\ t\in \R, &\label[3.1]\cr
        \mu(-t) &=-\mu(t) \quad \hbox{for all}\ t\in \R, &\label[3.2]\cr
        \mu(t)&=1 \quad \hbox{for}\ t\geq 1,  &\label[3.3]\cr
        \mu(t)&=-1 \quad \hbox{for}\ t\leq -1, &\label[3.4]\cr
        \mu'(t) &>0 \quad \hbox{for}\ t\in (-1,1). &\label[3.5]\cr}
    $$
It follows from \ref[3.1]--\ref[3.5] that
    $$  \eqalignno{
        &a_+(t), \ a_-(t) \in [1,3] \quad \hbox{for all}\ t\in\R, &\label[3.6]\cr
        &a_+(-t) = a_-(t),\ a_-(-t)=a_+(t)  \quad \hbox{for all}\ t\in\R, &\label[3.7]\cr
        &a_+'(t)=\mu'(t)>0, \ a_-'(t)=-\mu'(t) <0 \quad \hbox{for all}\ t\in(-1,1). &\label[3.8]\cr}
    $$
Finally we define $\varphi(t)\in C^1(\R,\R)$ by
    $$  \varphi(t)=\cases{  (t-1)^2 &for $t> 1$,\cr
                            0       &for $t\in [-1,1]$,\cr
                            (t+1)^2 &for $t<-1$.\cr}
    $$
We can see that $I(t,x_1,x_2,\cdots)$ has the following properties.

\medskip

\claim Proposition 3.1.
\item{(i)} $I(t,x_1,x_2,\cdots)\in C^1(X,\R)$;
\item{(ii)} $I(t,x_1,x_2,\cdots)$ is bounded from below and coercive on $X$;
\item{(iii)} $I(t,x_1,x_2,\cdots)$ satisfies $(PS)_c$ for all $c\in\R$;
\item{(iv)} $I(t,x_1,x_2,\cdots)$ is even in  $(t,x_1,x_2,\cdots)$, that is, 
$I(-t,-x_1,-x_2,\cdots)=I(t,x_1,x_2,\cdots)$;
\item{(v)} $I(t,x_1,x_2,\cdots)$ satisfies (A3').

\medskip

\claim Proof.
It follows from H\"older inequality that
    $$  \left|\sum_{j=1}^\infty b_j(x_j)^{3/2}\right| \leq
        \left(\sum_{j=1}^\infty b_j^4\right)^{1/4}\left(\sum_{j=1}^\infty \abs{x_j}^2\right)^{3/4},
                                                                        \eqno\label[3.9]
    $$
from which we can see that $I(t,x_1,x_2,\cdots)$ is well-defined as a functional on $X$.  Using
\ref[3.9], we can also see (i)--(iii).  
(iv) follows from \ref[3.7].

For $k\in\N$, setting $X^k=\{ (0,x_1,x_2,\cdots , x_k, 0, 0, \cdots);\, (x_1, x_2, \cdots, x_k) \in\R^k \}$, 
we can easily find $\rho_k>0$ such that
    $$  \sup\{ I(u);\, u\in X^k, \, \norm u_X=\rho_k\} <0.
    $$
Thus (v) holds.  \QED

\medskip

\noindent
By Proposition 3.1, we can apply the Clark Theorem to $I(t,x_1,x_2,\cdots)$.  On the other hand,
we have

\medskip

\proclaim Proposition 3.2.
Let $K$ be the set of all critical points of $I(t,x_1,x_2,\cdots)$.  Then we have
    $$  K = Z \cup N\cup (-N),                      \eqno\label[3.10]
    $$
where
    $$  \eqalignno{
        Z &= \{ (t,0,0,\cdots);\, t\in [-1,1]\},\cr
        N &= \{ (1, x_1, x_2, \cdots);\, x_j\in \{ 0, \, 9\cdot 3^{-2j},\, -3^{-2j}\} \ \hbox{for all}\ j\}. &\label[3.11]\cr}
    $$
Moreover we have
    $$  \eqalignno{
        I(u) &=0 \quad \hbox{for all}\ u\in Z,              &\label[3.12]\cr
        I(u) &<0 \quad \hbox{for all}\ u\in K\setminus Z.   &\label[3.13]\cr}
    $$

\medskip

\claim Proof.
Since 
    $$  \partial_t I(t,x_1,x_2,\cdots)=\varphi'(t) \quad \hbox{for all}\ (t,x_1,x_2,\cdots)\in X \ 
        \hbox{with} \ \abs t\geq 1,                             \eqno\label[3.14]
    $$
we see that $I'(t,x_1,x_2,\cdots)=0$ implies $t\in [-1,1]$.

In what follows, we assume that $(t,x_1,x_2,\cdots)$ ($t\in [-1,1]$) is a critical point of $I$ 
and give more precise description.  First we show

\smallskip

\noindent
{\sl Step 1:  For any $j\in\N$, 
    $$  x_j=0, \ 3^{-2j} a_+(t)^2, \ \hbox{or}\  -3^{-2j} a_-(t)^2. \eqno\label[3.15]
    $$
In particular, we have
\item{(i)} if $x_j\not=0$,
    $$  3^{-2j} \leq \abs{x_j} \leq 9\cdot 3^{-2j}.             \eqno\label[3.16]
    $$
\item{(ii)} $K\cap \{ t=1\} = N$,  $K\cap \{ t=-1\} = -N$, where $N$ is given in \ref[3.11].

}

\smallskip

\noindent
In fact, it follows from $\partial_{x_j}I(t,x_1,x_2,\cdots)=0$ that
    $$  x_j - 3^{-j}\left(a_+(t)(x_j)_+^{1/2}-a_-(t)(x_j)_-^{1/2}\right) =0.
    $$
From which we can get \ref[3.15].  (i) and (ii) follow from the property \ref[3.6] and
$a_+(1)=a_-(-1)=3$, $a_+(-1)=a_-(1)=1$.

\smallskip

\noindent
{\sl Step 2:  When $t\in (-1,1)$, it holds that $x_j=0$ for all $j$.}

\smallskip

\noindent
It follows from $\partial_t I(t,x_1,x_2,\cdots)=0$ that
    $$  \sum_{j=1}^\infty 3^{-j}\left( a_+'(t)(x_j)_+^{3/2} +a_-'(t)(x_j)_-^{3/2}\right) =0.
    $$
By \ref[3.8],
    $$  \sum_{j=1}^\infty 3^{-j}\left( (x_j)_+^{3/2} - (x_j)_-^{3/2}\right) =0.
    $$
Arguing indirectly, we assume that $x_j\not=0$ for some $j$ and let $j_0$ be the smallest integer
such that $x_j\not=0$.  Then we have
    $$  \eqalignno{
        3^{-j_0} \abs{x_{j_0}}^{3/2} &= \left| 3^{-j_0} \left( (x_{j_0})_+^{3/2} - (x_{j_0})_-^{3/2}\right) \right|
        =\left| \sum_{j=j_0+1}^\infty 3^{-j}\left( (x_j)_+^{3/2} - (x_j)_-^{3/2}\right)\right| \cr
        &\leq \sum_{j=j_0+1}^\infty 3^{-j} \abs{x_j}^{3/2}.             &\label[3.17]\cr}
    $$
By \ref[3.16],
    $$  \hbox{the right hand side of \ref[3.17]} 
        \leq \sum_{j=j_0+1}^\infty 3^{-j}(9\cdot 3^{-2j})^{3/2} = \sum_{j=j_0+1}^\infty 27\cdot 3^{-4j}
        < {2\over 3}\cdot 3^{-4j_0}.
    $$
Again by \ref[3.16], 
    $$  \hbox{the left hand side of \ref[3.17]} \geq 3^{-j_0} (3^{-2j_0})^{3/2} = 3^{-4j_0},
    $$
which is in contradiction with \ref[3.17].  Thus we have $x_j=0$ for all $j\in \N$.

\smallskip

{\sl
\noindent Step 3: Conclusion.}

\smallskip

\noindent \ref[3.10] follows from Steps 1--2.  We can also verify \ref[3.12]--\ref[3.13] easily.  \QED

\medskip

\noindent
As an immediate corollary to Proposition 3.2, we have

\medskip

\proclaim Corollary 3.3.
\item{(i)}  Points $(t,0,\cdots)\in X$ ($t\in (-1,1)$) cannot be accumulation points of critical
points with negative critical values.  
\item{(ii)} $(1,0,0,\cdots)$, $(-1,0,0,\cdots)\in X$ are accumulation points of critical points with
negative critical values.

\medskip

\noindent
Thus Corollary 3.3 shows that (C1) does not hold in general under the conditions (A1), (A2''), (A3').

\medskip

\BSS{3.2. An example which shows Theorem 1.8}
Next we give another example, which shows that without $(PS)_0$ the conclusion of Theorem 1.4 does not hold in
general.  
Here we work in the Hilbert space $(E,\norm\cdot_E)$ given by 
    $$  \eqalign{
        &E=H_0^1(0,1),\cr
        &\norm u_E=\left(\int_0^1\abs{u_x}^2\, dx\right)^{1/2} \quad \hbox{for}\ u\in E.\cr}
    $$
For $p\in (0,1)$ we define $J(u)\in C^1(E,\R)$ by
    $$  J(u)=\half\norm u_E^2-{1\over p+1}\int_0^1\abs u^{p+1}\, dx:\, E\to \R.
    $$
Critical points of $J(u)$ are solutions of the following sublinear elliptic equation:
    $$  \left\{ \eqalign{
        &u_{xx} + \abs u^{p-1} u =0 \quad \hbox{in}\ (0,1),\cr
        &u(0)=u(1)=0 \cr}
        \right.
    $$
and it has the following properties:
\item{(i)} $J(0)=0$, $J(u)$ is even in $u$,  bounded from below and coercive;
\item{(ii)} For any $k\in \N$, there exists a compact subset $A\subset E\setminus \{0\}$, which is
symmetric with respect to $0$, such that
    $$  \gamma(A)\geq k, \quad \max_{u\in A} J(u) <0.
    $$
Actually, for any $k$-dimensional subspace $H\subset E$,
    $$  A=\{ u\in H;\, \norm u_E=\rho\}
    $$
with small $\rho>0$ gives the desired compact set.
\item{(iii)} $J(u)$ satisfies $(PS)_c$ for all $c\in\R$.

\medskip

\noindent
We define $I(u):\, E\to \R$ by
    $$  I(u)=\cases{    1-\cos(2\pi \norm u_E^2)    &for $\norm u_E\leq 1$,\cr
                        J((\norm u_E^2-1)u)         &for $\norm u_E>1$,\cr}
    $$
Then we have

\medskip

\proclaim Proposition 3.4.
$I(u):\, E\to \R$ satisfies the assumptions of Theorem 1.1.  However $0$ is an isolated critical point
and the conclusion of Theorem 1.4 does not hold.

\medskip

\claim Proof.
Clearly $I(u)$ is even, bounded from below and coercive.  Moreover $I(u)$ also satisfies $(PS)_c$ for all $c<0$.
In fact, if $(u_j)_{j=1}^\infty$ satisfies $I(u_j)\to c<0$ and $\norm{I'(u_j)}\to 0$, then we can easily see
that $(u_j)_{j=1}^\infty$ is bounded as $j\to\infty$ and, after taking a subsequence, we may assume that
$\norm{u_j}_E\to d$ for some $d>1$.  Using this fact, we can see that  $(u_j)_{j=1}^\infty$ has a strongly
convergent subsequence.
(Since all points on the unit sphere $S=\{ x\in E;\ \norm u_E=1\}$ are critical points 
of $I(u)$ with critical value $0$ and $S$ is not compact, we note that $(PS)_0$ fails.)
Thus $I(u)$ satisfies (A1) and (A2).  We can see that (A3) holds easily.  
In fact, for any $k$-dimensional subspace $H\subset E$, choosing $\rho>0$ small,
    $$  A=\{ u\in E;\, \norm u_E= 1+\rho\}
    $$
satisfies $\gamma(A)\geq k$ and $\sup_{u\in A}I(u)<0$.  Thus $I(u)$ satisfies the assumptions
of Theorem 1.1.

Clearly $0$ is an isolated critical point of $I(u)$ and $I(u)$ does not have a sequence of critical points 
with \ref[1.2]--\ref[1.4].  \QED

\bigskip

\noindent
{\bf Acknowledgments.}
This research is motivated by Professor Zhaoli Liu's lecture in Capital 
Normal University, to which the first author attended.
Authors would like to thank Professor Liu for suggesting them
to study such an interesting problem.
Authors would also like to thank Professor Ryuji Kajikiya for helpful discussions.

Authors started this research during the first and the third authors' visit to Department of Mathematics,
School of Science and Engineering, Waseda University. They would like to gratefully acknowledge 
Waseda University for cordial invitations and hospitality. The first author wishes to express his 
gratitude to the support of Graduate School of Capital Normal University. The third author is greatly
indebted to China Scholarship Council for their support.

The second author is partially supported by JSPS Grants-in-Aid for Scientific Research (B) (25287025) and 
Waseda University Grant for Special Research Projects 2016B-120. 

\bigskip


\centerline{\bf References}

\medskip

{\parindent=1.5\parindent

\bibitem[Cl] D. C. Clark, 
A variant of the Lusternik-Schnirelman theory, Indiana Univ. Math. J. 22 (1972), 65--74.

\bibitem[Co] C. V. Coffman, 
A minimum-maximum principle for a class of non-linear integral equations, 
J. Analyse Math. 22 (1969), 391--419 .



\bibitem[H] H.-P. Heinz, 
Free Ljusternik-Schnirelman theory and the bifurcation diagrams of certain singular nonlinear problems, 
J. Differential Equations 66 (1987), no. 2, 263--300.

\bibitem[Ka1] R. Kajikiya, 
A critical point theorem related to the symmetric mountain pass lemma and its applications 
to elliptic equations, J. Funct. Anal. 225 (2005), no. 2, 352--370.

\bibitem[Ka2] R. Kajikiya,
private communication, 2017.

\bibitem[Kr] M. A. Krasnosel'skii, 
Topological methods in the theory of nonlinear integral equations. 
Translated by A. H. Armstrong; translation edited by J. Burlak. 
A Pergamon Press Book The Macmillan Co., New York 1964 xi + 395 pp. 

\bibitem[LW] Z. Liu, Z.-Q. Wang, 
On Clark's theorem and its applications to partially sublinear problems, 
Ann. Inst. H. Poincar\'e Anal. Non Lin\'eaire 32 (2015), no. 5, 1015--1037.

\bibitem[R] P. H. Rabinowitz, 
Minimax methods in critical point theory with applications to differential equations. 
CBMS Regional Conference Series in Mathematics, 65. 
Published for the Conference Board of the Mathematical Sciences, Washington, DC; 
by the American Mathematical Society, Providence, RI, 1986. viii+100 pp.

}

\bye